\documentclass[10pt]{amsart}

\usepackage{amsmath,amsfonts,amssymb,amsthm,amscd,epsfig}

\def\C{\mathbb{C}}
\def\Z{\mathbb{Z}}

\def\Q{\mathbb{Q}}

\def\g{\ensuremath{\mathfrak{g}}}

\def\V{\mathbf{V}}
\def\W{\mathbf{W}}
\def\D{\mathbf{D}}
\def\v{\mathbf{v}}

\def\d{\mathbf{d}}
\def\e{\mathbf{e}}

\def\L{\mathfrak{L}}

\def\d{{\mathbf{d}}}

\def\ke{{\tilde e}}
\def\kf{{\tilde f}}

\newcommand{\comment}[1]{}

\DeclareMathOperator{\im}{Im} %Image of a map
\DeclareMathOperator{\Hom}{Hom}

\DeclareMathOperator{\End}{End}
\DeclareMathOperator{\Aut}{Aut}
\DeclareMathOperator{\inc}{in}
\DeclareMathOperator{\out}{out}

\DeclareMathOperator{\tr}{tr}
\DeclareMathOperator{\Irr}{Irr}

\DeclareMathOperator{\wt}{wt}

\newtheorem{theo}{Theorem}[section]
\newtheorem{prop}[theo]{Proposition}
\newtheorem{lem}[theo]{Lemma}
\newtheorem{cor}[theo]{Corollary}

\newtheorem{defin}[theo]{Definition}
\newtheorem*{rem*}{Remark}

\numberwithin{equation}{section}

\begin{document}
\title{Quiver Varieties and Demazure modules}
\author{Alistair Savage}
\address{The Fields Institute for Research in Mathematical Sciences
and The University of Toronto \\
Toronto, Ontario \\ Canada} \email{alistair.savage@aya.yale.edu}
\thanks{This research was supported in part by the Natural
Sciences and Engineering Research Council (NSERC) of Canada}
\subjclass[2000]{Primary 16G20,17B37}
%\date{June 29, 2005}

\begin{abstract}
Using subvarieties, which we call Demazure quiver varieties, of
the quiver varieties of Nakajima, we give a geometric realization
of Demazure modules of Kac-Moody algebras with symmetric Cartan
data. We give a natural geometric characterization of the extremal
weights of a representation and show that Lusztig's semicanonical
basis is compatible with the filtration of a representation by
Demazure modules.  For the case of $\widehat{\mathfrak{sl}}_2$, we
give a characterization of the Demazure quiver variety in terms of
a nilpotency condition on quiver representations and an explicit
combinatorial description of the Demazure crystal in terms of
Young pyramids.
\end{abstract}

\maketitle

\section{Introduction}
\label{sec:dem-def}

Let $\g$ be a symmetrizable Kac-Moody Lie algebra.  Let $W$ be its
Weyl group with Bruhat order $\prec$ and $\{r_i\}_{i \in I}$ the
set of simple reflections. Let $U(\g)$ be the universal enveloping
algebra of $\g$ and $U^+(\g)$ the subalgebra generated by the
$e_i$'s. For a dominant integral weight $\lambda$, let
$V(\lambda)$ be the irreducible highest weight $\g$-module with
highest weight $\lambda$.  For $w \in W$, it is known that the
extremal weight space $V(\lambda)_{w \lambda}$ is one dimensional.
Here $V_\mu$ denotes the $\mu$-weight space of a $\g$-module $V$.
Let $V_w(\lambda)$ denote the $U^+(\g)$-module generated by
$V(\lambda)_{w \lambda}$.  The modules $V_w(\lambda)$, $w \in W$,
are called \emph{Demazure modules}.  They are finite dimensional
subspaces which form a filtration of $V(\lambda)$ that is
compatible with the Bruhat order of $W$.  That is, $V_w(\lambda)
\subset V_{w'}(\lambda)$ whenever $w \prec w'$ and $\bigcup_{w \in
W} V_w(\lambda) = V(\lambda)$.

An interesting property of Demazure modules is their relation to
integrable models in statistical mechanics. The so-called
``one-point functions" in exactly solvable two-dimensional
lattices can be evaluated using Baxter's corner transfer matrix
method.  This method reduces the computation of one-point
functions to a weighted sum over combinatorial objects called
paths which are one-dimensional configurations defined on the half
line.  The generating function of all paths turns out to be the
character of a highest weight module of an affine Lie algebra. The
same physical objects have been evaluated in a completely
different way, using Bethe Ansatz methods.  Equating the two
expressions for the same object, one obtains Rogers-Ramanujan-type
identities between $q$-series.  When attempting to provide
rigorous proofs of these identities, it is often necessary to work
at the level of the half-line of finite length $L$ and then take
the limit $L \to \infty$.  A natural question then arises as to
the representation theoretic meaning of this finite length.  In
\cite{FMO98}, it was shown that the $L$-restricted generating
functions are closely related to the characters of Demazure
modules.

Another interesting property of Demazure modules is their
compatibility with the global basis and crystal structure.  Let
$V^q(\lambda)$ and $V^q_w(\lambda)$ be the $q$-analogues of
$V(\lambda)$ and $V_w(\lambda)$ respectively (i.e. simply replace
all objects in the above definition with their $q$-analogues). Let
$(L(\lambda), B(\lambda))$ be the crystal base of $V^q(\lambda)$.
Kashiwara showed in \cite{Kas93} that for each $w \in W$, there
exists a subset $B_w(\lambda)$ of $B(\lambda)$ such that
\begin{equation}
\label{eq:demazure-crystalbase}
\frac{V^q_w(\lambda) \cap L(\lambda)}{V^q_w(\lambda) \cap
qL(\lambda)} = \sum_{b \in B_w(\lambda)} \Q b.
\end{equation}
Kashiwara also proves that
\[
V^q_w(\lambda) = \bigoplus_{b \in B_w(\lambda)} \Q(q)
G_\lambda(b),
\]
where $\{G_\lambda(b)\}_{b \in B(\lambda)}$ is the lower global
basis.  In addition, $B_w(\lambda)$ has the following recursive
property: If $r_i w \succ w$, then $B_{r_i w} (\lambda) =
\bigcup_{n \ge 0} \kf_i^n B_w(\lambda) \backslash \{0\}$.  In
particular, if $w = r_{i_l} \dots r_{i_1}$ is a reduced expression
of $w$, then
\begin{equation}
\label{eq:dqv}
B_w(\lambda) = \{\kf_{i_l}^{n_l} \dots \kf_{i_1}^{n_1} b_\lambda\
|\ n_1,\dots,n_l \in \Z_{\ge 0}\} \backslash \{0\}.
\end{equation}

Lusztig \cite{L91} and Nakajima \cite{N94,N98} have developed a
geometric realization of universal enveloping algebras of
Kac-Moody algebras with symmetric Cartan data and their integrable
highest-weight representations.  The algebraic objects are
realized as a set of perverse sheaves or constructible functions
on varieties attached to quivers or in the homology of these
varieties.  One actually considers a union of many quiver
varieties, each one corresponding to a given weight space of the
algebra or representation.  The crystal structure of the algebras
and their representations has also been realized on the set of
irreducible components of these varieties (cf. \cite{KS97,S02}).
It is known (cf. \cite{GL93}) that the lower global basis defined
by Kashiwara coincides with Lusztig's canonical basis, defined
geometrically by considering perverse sheaves on quiver varieties.
Given these developments, it is natural to ask for an
interpretation of Demazure modules in terms of quiver varieties.
This is the theme of the current paper.

We find that the Demazure modules are very natural in this
geometric setting.  In particular, it turns out that the quiver
varieties corresponding to a given representation of highest
weight $\lambda$ are single points exactly for the varieties
corresponding to the extremal weights $w \lambda$, $w \in W$.
These points, which we denote $[x_w,t_w]$, are orbits of
representations of the quiver.  If one then restricts the quiver
varieties by only considering subrepresentations of the
representative $(x_w,t_w)$ of $[x_w,t_w]$ (we call the resulting
variety the \emph{Demazure quiver variety}), we obtain a geometric
realization of the Demazure modules.

In the geometric realization of representations using
constructible functions, there is a distinguished function
attached to each irreducible component of the quiver varieties.
These functions form a special basis of the representation, called
the semicanonical basis.  This basis has many nice properties such
as compatibility with various filtrations and with canonical
antiautomorphisms of the universal enveloping algebras (cf.
\cite{Lus00}). In this paper, we add another property to this
list. Namely, we show that the semicanonical basis is compatible
with the filtration by Demazure modules.  More specifically, let
$B$ be the set of irreducible components of the quiver varieties
corresponding to a representation of highest weight $\lambda$ and
$B_w$ the subset corresponding to the Demazure crystal for some $w
\in W$ (we show that this is exactly the set of irreducible
components completely contained in the Demazure quiver variety).
We then show that
\[
V_w(\lambda) \cong \bigoplus_{X \in B_w} \C g_X,
\]
where the $g_X$ is the semicanonical basis element corresponding
to the irreducible component $X$.

Finally, we consider the case of $\g = \widehat{\mathfrak{sl}}_2$
in more detail.  We show that in this case, the Demazure quiver
variety is simply given by imposing nilpotency of a certain order
on the quiver representations.  We also give an explicit
description of the Demazure crystal in terms of the Young pyramids
defined in \cite{Sav03b}.  These objects were inspired by the
geometric realization of crystals.

The organization of the paper is as follows.  In
Sections~\ref{sec:lus_def}, \ref{sec:def_nak} and
\ref{sec:structure} we review the geometric construction of
representations of Kac-Moody algebras using quiver varieties.  In
Section~\ref{sec:extremal} we give a geometric characterization of
the extremal weights of these representations.  We define the
Demazure quiver variety and examine the geometric Demazure crystal
in Section~\ref{sec:dqv} and then prove the compatibility of the
semicanonical basis with the Demazure filtration in
Section~\ref{sec:sc-basis}.  Finally, in
Sections~\ref{sec:pyramids} and \ref{sec:sl2-dqv} we consider the
special case of $\widehat{\mathfrak{sl}}_2$, describing the
Demazure crystals in terms of Young pyramids and the Demazure
quiver varieties in terms of the nilpotency of quiver
representations respectively.

The author would like to thank M. Okado for raising the question
of a quiver variety interpretation of Demazure modules, Kyoto and
Osaka universities for their hospitality and T. Miwa for the
invitation to visit these institutions.  He would also like to
thank I. B. Frenkel for useful discussions and suggestions and
note that in \cite{FKS03} I. B. Frenkel, M. Khovanov and O.
Schiffmann hypothesized the existence of a simple geometric
description of Demazure modules using quiver varieties.

%%%%%%%%%%%%%%%%%%%%%%%%%%%%%%%%%%%%%%%%%%%%%%%%%%%%%%%%%%%%%%%%%%%%%%%

\section{Lusztig's quiver variety}
\label{sec:lus_def}

In this section, we will recount the description given in
\cite{L91} of Lusztig's quiver variety. See this reference for
details, including proofs.

Let $I$ be the set of vertices of the Dynkin graph of a symmetric
Kac-Moody Lie algebra $\mathfrak{g}$ and let $H$ be the set of
pairs consisting of an edge together with an orientation of it.
For $h \in H$, let $\inc(h)$ (resp. $\out(h)$) be the incoming
(resp. outgoing) vertex of $h$.  We define the involution $\bar{\
}: H \to H$ to be the function which takes $h \in H$ to the
element of $H$ consisting of the same edge with opposite
orientation.  An \emph{orientation} of our graph is a choice of a
subset $\Omega \subset H$ such that $\Omega \cup \bar{\Omega} = H$
and $\Omega \cap \bar{\Omega} = \emptyset$.

Let $\mathcal{V}$ be the category of finite-dimensional $I$-graded
vector spaces $\V = \oplus_{i
  \in I} \V_i$ over $\C$ with morphisms being linear maps
respecting the grading.  Then $\V \in \mathcal{V}$ shall denote
that $\V$ is an object of $\mathcal{V}$.  We associate to the
graded dimension $\v = (\v_i)_{i \in I}$ of $\V$ the element
$\alpha_\v = \sum_{i \in I} \v_i \alpha_i$ of the root lattice of
$\mathfrak{g}$. Here the $\alpha_i$ are the simple roots
corresponding to the vertices of our quiver (graph with
orientation), whose underlying graph is the Dynkin graph of
$\mathfrak{g}$.

Given $\V \in \mathcal{V}$, let
\[
\mathbf{E_V} = \bigoplus_{h \in H} \Hom (\V_{\out(h)},
\V_{\inc(h)}).
\]
For any subset $H' \subset H$, let $\mathbf{E}_{\V, H'}$ be the
subspace of $\mathbf{E_V}$ consisting of all vectors $x = (x_h)$
such that $x_h=0$ whenever $h \not\in H'$.  The algebraic group
$G_\V = \prod_i \Aut(\V_i)$ acts on $\mathbf{E_V}$ and
$\mathbf{E}_{\V, H'}$ by
\[
(g,x) = ((g_i), (x_h)) \mapsto (g_{\inc(h)} x_h g_{\out(h)}^{-1}).
\]

Define the function $\varepsilon : H \to \{-1,1\}$ by $\varepsilon
(h) = 1$ for all $h \in \Omega$ and $\varepsilon(h) = -1$ for all
$h \in {\bar{\Omega}}$.  The Lie algebra of $G_\V$ is
$\mathbf{gl_V} = \prod_i \End(\V_i)$ and it acts on $\mathbf{E_V}$
by
\[
(a,x) = ((a_i), (x_h)) \mapsto [a,x] = (x'_h) = (a_{\inc(h)}x_h -
x_h a_{\out(h)}).
\]
Let $\left<\cdot,\cdot\right>$ be the nondegenerate,
$G_\V$-invariant, symplectic form on $\mathbf{E_V}$ with values in
$\C$ defined by
\[
\left<x,y\right> = \sum_{h \in H} \varepsilon(h) \tr (x_h
y_{\bar{h}}).
\]
Note that $\mathbf{E_V}$ can be considered as the cotangent space
of $\mathbf{E}_{\V, \Omega}$ under this form.

The moment map associated to the $G_{\V}$-action on the symplectic
vector space $\mathbf{E_V}$ is the map $\psi : \mathbf{E_V} \to
\mathbf{gl_V}$ with $i$-component $\psi_i : \mathbf{E_V} \to \End
\V_i$ given by
\[
\psi_i(x) = \sum_{h \in H,\, \inc(h)=i} \varepsilon(h) x_h
x_{\bar{h}} .
\]

\begin{defin}[\cite{L91}]
\label{def:nilpotent}
An element $x \in \mathbf{E_V}$ is said to be \emph{nilpotent} if
there exists an $N \ge 1$ such that for any sequence $h_1, h_2,
\dots, h_N$ in $H$ satisfying $\inc (h_1) = \out (h_2)$, $\inc
(h_2) = \out (h_3)$, \dots, $\inc (h_{N-1}) = \out (h_N)$, the
composition $x_{h_N} \dots x_{h_2} x_{h_1} : \V_{\out (h_1)} \to
  \V_{\inc (h_N)}$ is zero.
\end{defin}

\begin{defin}[\cite{L91}] Let $\Lambda_\V$ be the set of all
  nilpotent elements $x \in \mathbf{E_V}$ such that $\psi_i(x) = 0$
  for all $i \in I$.  We call $\Lambda_\V$ \emph{Lusztig's quiver
  variety}.
\end{defin}

It is known that the irreducible components of $\Lambda_\V$ are in
one-to-one correspondence with a basis of the $\alpha_\v$-weight
space of $U^+(\g)$.

%%%%%%%%%%%%%%%%%%%%%%%%%%%%%%%%%%%%%%%%%%%%%%%%%%%%%%%%%%%%%%%%%%%%

\section{Nakajima's quiver variety}
\label{sec:def_nak}

We introduce here a description of the quiver varieties first
presented in \cite{N94}.  See \cite{N94} and \cite{N98} for
details.

\begin{defin}[\cite{N94}]
\label{def:lambda} For $\v, \d \in (\Z_{\ge 0})^I$, choose
$I$-graded vector spaces $\V$ and $\D$ of graded dimensions $\v$
and $\d$ respectively.  We associate to $\d = (\d_i)_{i \in I}$
the element $\lambda_\d = \sum_i \d_i \omega_i$ of the root
lattice of \g, where the $\omega_i$ are the fundamental weights of
\g. Recall that we associated to $\v$ the weight $\alpha_\v =
\sum_i \v_i \alpha_i$. Then define
\[
\Lambda \equiv \Lambda(\v,\d) = \Lambda_\V \times \sum_{i \in I}
\Hom (\V_i, \D_i).
\]
\end{defin}

Now, suppose that $\mathbf{S}$ is an $I$-graded subspace of $\V$.
For $x \in \Lambda_\V$ we say that $\mathbf{S}$ is
\emph{$x$-stable} if $x(\mathbf{S}) \subset \mathbf{S}$.

\begin{defin}[\cite{N94}]
\label{def:lambda-stable} Let $\Lambda^{\text{st}} =
\Lambda(\v,\d)^{\text{st}}$ be the set of all $(x, t) \in
\Lambda(\v,\d)$ satisfying the following condition:  If
$\mathbf{S}=(\mathbf{S}_i)$ with $\mathbf{S}_i \subset \V_i$ is
$x$-stable and $t_i(\mathbf{S}_i) = 0$ for $i \in I$, then
$\mathbf{S}_i = 0$ for $i \in I$.
\end{defin}

The group $G_\V$ acts on $\Lambda(\v,\d)$ via
\[
(g,(x,t)) = ((g_i), ((x_h), (t_i))) \mapsto ((g_{\inc (h)} x_h
g_{\out (h)}^{-1}), (t_i g_i^{-1})).
\]
and the stabilizer of any point of $\Lambda(\v,\d)^{\text{st}}$ in
$G_{\V}$ is trivial
  (see \cite[Lemma~3.10]{N98}).  We then make the following definition.
\begin{defin}[\cite{N94}]
\label{def:L} Let $\L \equiv \L(\v,\d) = \Lambda(\v,\d)^{\text{st}}
/ G_{\V}$.  We call this \emph{Nakajima's quiver variety}.
\end{defin}

\begin{lem}[\cite{N98}]
\label{lem:nak-dim}
We have
\[
\dim_\C \L(\v,\d) = \frac{1}{2} \v^t (2\d - C\v),
\]
where $C$ is the Cartan matrix of $\g$.
\end{lem}

%%%%%%%%%%%%%%%%%%%%%%%%%%%%%%%%%%%%%%%%%%%%%%%%%%%%%%%%%%%%%%%%%%%%

\section{The Lie Algebra Action}
\label{sec:structure}

We summarize here some results from \cite{N94} that will be needed
in the sequel.  See this reference for more details, including
proofs. We keep the notation of Sections~\ref{sec:lus_def}
and~\ref{sec:def_nak}.

Let $\mathbf{d, v, v', v''} \in \Z_{\ge 0}^I$ be such that
$\mathbf{v} = \mathbf{v'} + \mathbf{v''}$.  Consider the maps
\begin{equation}
\label{eq:diag_action}
\Lambda(\mathbf{v}'',\mathbf{0}) \times
\Lambda(\mathbf{v}',\mathbf{d}) \stackrel{p_1}{\leftarrow}
\mathbf{\tilde F (v,d;v'')} \stackrel{p_2}{\rightarrow}
\mathbf{F(v,d;v'')} \stackrel{p_3}{\rightarrow}
\Lambda(\mathbf{v},\mathbf{d}),
\end{equation}
where the notation is as follows.  A point of
$\mathbf{F(v,d;v'')}$ is a point $(x,t) \in
\Lambda(\mathbf{v},\mathbf{d})$ together with an $I$-graded,
$x$-stable subspace $\mathbf{S}$ of $\mathbf{V}$ such that $\dim
\mathbf{S} = \mathbf{v'} = \mathbf{v} - \mathbf{v''}$.  A point of
$\mathbf{\tilde
  F (v,d;v'')}$ is a point $(x,t,\mathbf{S})$ of $\mathbf{F(v,d;v'')}$
together with a collection of isomorphisms $R'_i : \mathbf{V}'_i
\cong \mathbf{S}_i$ and $R''_i : \mathbf{V}''_i \cong \mathbf{V}_i
/ \mathbf{S}_i$ for each $i \in I$.  Then we define
$p_2(x,t,\mathbf{S}, R',R'') = (x,t,\mathbf{S})$,
$p_3(x,t,\mathbf{S}) = (x,t)$ and $p_1(x,t,\mathbf{S},R',R'') =
(x'',x',t')$ where $x'', x', t'$ are determined by
\begin{align*}
R'_{\inc(h)} x'_h &= x_h R'_{\out(h)} : \mathbf{V}'_{\out(h)} \to
\mathbf{S}_{\inc(h)}, \\
t'_i &= t_i R'_i : \mathbf{V}'_i \to \D_i \\
R''_{\inc(h)} x''_h &= x_h R''_{\out(h)} : \mathbf{V}''_{\out(h)}
\to \mathbf{V}_{\inc(h)} / \mathbf{S}_{\inc(h)}.
\end{align*}
It follows that $x'$ and $x''$ are nilpotent.

\begin{lem}[{\cite[Lemma 10.3]{N94}}]
One has
\[
(p_3 \circ p_2)^{-1} (\Lambda(\mathbf{v},\mathbf{d})^{\text{st}})
\subset p_1^{-1} (\Lambda(\mathbf{v}'',\mathbf{0}) \times
\Lambda(\mathbf{v}',\mathbf{d})^{\text{st}}).
\]
\end{lem}

Thus, we can restrict \eqref{eq:diag_action} to
$\Lambda^{\text{st}}$, forget the
$\Lambda(\mathbf{v}'',\mathbf{0})$-factor and consider the
quotient by $G_\mathbf{V}$, $G_\mathbf{V'}$.  This yields the
diagram
\begin{equation}
\label{eq:diag_action_mod} \L(\mathbf{v'}, \mathbf{d})
\stackrel{\pi_1}{\leftarrow} \mathfrak{F}(\mathbf{v}, \mathbf{d};
\mathbf{v} - \mathbf{v'}) \stackrel{\pi_2}{\rightarrow}
\L(\mathbf{v}, \mathbf{d}),
\end{equation}
where
\[
\mathfrak{F}(\mathbf{v}, \mathbf{d}, \mathbf{v} - \mathbf{v'})
\stackrel{\text{def}}{=} \{ (x,t,\mathbf{S}) \in
\mathbf{F(v,d;v-v')}\,
  |\, (x,t) \in \Lambda(\mathbf{v},\mathbf{d})^{\text{st}} \} / G_\mathbf{V}.
\]

Let $M(\L(\mathbf{v}, \mathbf{d}))$ be the vector space of all
constructible functions on $\L(\mathbf{v}, \mathbf{d})$. For a
subvariety $Y$ of a variety $A$, let $\mathbf{1}_Y$ denote the
function on $A$ which takes the value 1 on $Y$ and 0 elsewhere.  Let
$\chi (Y)$ denote the Euler characteristic of the algebraic variety
$Y$.  Then for a map $\pi$ between algebraic varieties $A$ and $B$,
let $\pi_!$ denote the map between the abelian groups of
constructible functions on $A$ and $B$ given by
\[
\pi_! (\mathbf{1}_Y)(y) = \chi (\pi^{-1}(y) \cap Y),\ Y \subset A,
\]
and let $\pi^*$ be the pullback map from functions on $B$ to
functions on $A$ acting as $\pi^* f(y) = f(\pi(y))$. Then define
\begin{align*}
&H_i : M(\L(\mathbf{v}, \mathbf{d})) \to
M(\L(\mathbf{v}, \mathbf{d})); \quad H_i f = u_i f, \\
&E_i : M(\L(\mathbf{v}, \mathbf{d})) \to M(\L(\mathbf{v} -
\mathbf{e}^i, \mathbf{d})); \quad E_i f =
(\pi_1)_! (\pi_2^* f), \\
&F_i : M(\L(\mathbf{v} - \mathbf{e}^i, \mathbf{d})) \to
M(\L(\mathbf{v}, \mathbf{d})); \quad F_i g = (\pi_2)_! (\pi_1^* g).
\end{align*}
Here
\[
\mathbf{u} = {^t(u_0, \dots, u_n)} = \mathbf{d} - C \mathbf{v}
\]
where $C$ is the Cartan matrix of $\mathfrak{g}$ and we are using
diagram~\eqref{eq:diag_action_mod} with $\mathbf{v}' = \mathbf{v}
- \mathbf{e}^i$ where $\mathbf{e}^i$ is the vector whose
components are given by $\mathbf{e}^i_j = \delta_{ij}$.

Now let $\varphi$ be the constant function on $\L(\mathbf{0},
\mathbf{d})$ with value 1.  Let $L(\mathbf{d})$ be the vector space
of functions generated by acting on $\varphi$ with all possible
combinations of the operators $F_i$.  Then let
$L(\mathbf{v},\mathbf{d}) = M(\L(\mathbf{v}, \mathbf{d})) \cap
L(\mathbf{d})$.

\begin{prop}[{\cite[Thm 10.14]{N94}}]
The operators $E_i$, $F_i$, $H_i$ on $L(\mathbf{d})$ provide the
structure of the irreducible highest weight integrable
representation of $\mathfrak{g}$ with highest weight $\lambda_\d$.
Each summand of the decomposition $L(\mathbf{d}) =
\bigoplus_\mathbf{v} L(\mathbf{v}, \mathbf{d})$ is a weight space
with weight $\lambda_\d - \alpha_\v$.
\end{prop}

Let $X \in \Irr \L(\mathbf{v}, \mathbf{d})$, the set of irreducible
components of $\L(\v,\d)$, and define a linear map $T_X :
L(\mathbf{v}, \mathbf{d}) \to \C$ as in \cite[\S 3.8]{L92}. The map
$T_X$ associates to a constructible function $f \in L(\mathbf{v},
\mathbf{d})$ the (constant) value of $f$ on a suitable open dense
subset of $X$. The fact that $L(\mathbf{v}, \mathbf{d})$ is
finite-dimensional allows us to take such an open set on which
\emph{any} $f \in L(\mathbf{v}, \mathbf{d})$ is constant.  So we
have a linear map
\[
\Phi : L(\mathbf{v}, \mathbf{d}) \to \C^{\Irr \L(\mathbf{v},
  \mathbf{d})}.
\]
The following proposition is proved in \cite[4.16]{L92} (slightly
generalized in \cite[Proposition 10.15]{N94}).

\begin{prop}
\label{prop:func_irrcomp_isom} The map $\Phi$ is an isomorphism; for
any $X \in \Irr \L(\mathbf{v}, \mathbf{d})$, there is a function
$g_X \in L(\mathbf{v}, \mathbf{d})$ such that for some open dense
subset $O$ of $X$ we have $g_X|_O = 1$ and for some closed
$G_\mathbf{V}$-invariant subset $K \subset \L(\mathbf{v},
\mathbf{d})$ of dimension $< \dim \L(\mathbf{v}, \mathbf{d})$ we
have $g_X=0$ outside $X \cup K$.  The functions $g_X$ for $X \in
\Irr \Lambda(\mathbf{v},\mathbf{d})$ form a basis of
$L(\mathbf{v},\mathbf{d})$.
\end{prop}

We note that in \cite{N94} it was only asserted that the map
$\Phi$ is an isomorphism for $\g$ of type $A$, $D$, $E$ or of
affine type.  The extension to arbitrary (symmetric) type follows
from the results of \cite{KS97,S02,Lus00}.  The basis given by the
$g_X$ is called the semicanonical basis and has some very nice
properties (see \cite{Lus00}).

The set of irreducible components of Nakajima's quiver variety can
be endowed with the structure of a crystal.  We denote this
crystal by $B(\d)$.  It is isomorphic to the crystal
$B(\lambda_\d)$ of the module $V(\lambda_\d)$. We refer the reader
to \cite{KS97,S02} for the details of this construction.

%%%%%%%%%%%%%%%%%%%%%%%%%%%%%%%%%%%%%%%%%%%%%%%%%%%%%%%%%%%%%%%%%%%%%%%

\section{Quiver varieties and extremal weights}
\label{sec:extremal}

We give here a natural geometric characterization of the extremal
weights $w \lambda$, $w \in W$, of the module $V(\lambda)$.

\begin{prop}
\label{prop:qv-point}
The quiver variety $\L(\v,\d)$ is a point if and only if
$\lambda_\d - \alpha_\v = w \lambda_\d$ for some $w \in W$.
\end{prop}

\begin{proof}
We first prove the ``if" part of the statement.  Since it follows
from the definition of Nakajima's quiver varieties that $L(0,\d)$
is a point, it suffices to show that if $\L(\v,\d)$ is a point
where $\lambda_\d - \alpha_\v = w \lambda_\d$ for some $w \in W$
then $\L(\v',\d)$ is a point where $\lambda_\d - \alpha_{\v'} =
r_i w \lambda_\d$ for an arbitrary simple reflection $r_i$.  Let
$\mu = \lambda_\d - \alpha_\v$.  Then
\[
r_i \mu = \mu - \left< \mu,\alpha_i \right> \alpha_i = \sum_j \d_j
\omega_j - \sum_j \v_j \alpha_j - \left<\mu, \alpha_i \right>
\alpha_i.
\]
So
\begin{align*}
\v_j' &= \v_j \text{ for } j \ne i, \\
\v'_i &= \v_i + \left< \mu, \alpha_i \right> \\
&= \v_i + \d_i - \sum_j \v_j \left< \alpha_j, \alpha_i \right> \\
&= \v_i + \d_i - \sum_j \v_j C_{ji} \\
&= \v_i + \d_i - \v^t C\e^i,
\end{align*}
where $C$ is the Cartan matrix of $\g$. Let $\delta = \v^t C\e^i$.
Note that since $C$ is symmetric, we also have $\delta = (\e^i)^t
C\v$. Then by Lemma~\ref{lem:nak-dim},
\begin{align*}
2\dim_\C \L(\v',\d) &= (\v')^t (2\d - C\v') \\
&= (\v + (\d_i - \delta)\e^i)^t (2\d - C(\v+(\d_i - \delta)\e^i)) \\
&= \v^t(2\d - C\v) - (\d_i - \delta)\v^tC\e^i + 2(\d_i -
\delta)(\e^i)^t \d \\
&\quad \quad - (\d_i-\delta)(\e^i)^tC\v - (\d_i-\delta)^2 (\e^i)^t
C\e^i.
\end{align*}
Since $\v^t(2\d - C\v) = 2\dim_\C \L(\v,\d) = 0$, we have
\begin{align*}
2\dim_\C \L(\v',\d) &= -2(\d_i - \delta)\delta + 2(\d_i -
\delta)\d_i
-2(\d_i - \delta)^2 \\
&=0.
\end{align*}
Since it is known that all $\L(\v,\d)$ are connected, it follows
that $\L(\v',\d)$ is a point.

Now, to prove the ``only if" part of the proposition, we use the
crystal structure on the set of irreducible components of
Nakajima's quiver variety.  Let $X$ be an irreducible component of
$\L(\v,\d)$ of dimension zero (i.e. $X$ is a point).  Since we
know that $\L(0,\d)$ is a point, we can assume $\v \ne 0$.  Choose
an $i \in I$ such that $\ke_i X \ne 0$ (since the crystal for
$V(\lambda_\d)$ is connected, we can always find such an $i$). Let
$c$ be the maximum exponent such that $\ke_i^c X \ne 0$. Let
$\ke_i^c X$ be an irreducible component of $\L(\v',\d)$. That is,
\begin{equation}
\label{eq:irrcomp-wd}
\lambda_\d - \alpha_{\v'} = \lambda_\d - \alpha_\v + c\alpha_i.
\end{equation}
It follows from \cite[Lemma~4.2.2]{S02} that $\dim \ke_i^c X \le
\dim X = 0$.  So $\dim \ke_i^c X = 0$.  Thus, it also follows from
\cite[Lemma~4.2.2]{S02} that
\[
c = \left< h_i, \lambda_\d - \alpha_{\v'} \right> =\left<
\alpha_i, \lambda_\d - \alpha_{\v'} \right>
\]
Substituting into \eqref{eq:irrcomp-wd}, we see that
\[
\lambda_\d - \alpha_\v = r_i (\lambda_\d - \alpha_{\v'}).
\]
The result then follows by induction on the length of $w$ or,
equivalently, on $|\v| = \sum_i \v_i$.
\end{proof}

%%%%%%%%%%%%%%%%%%%%%%%%%%%%%%%%%%%%%%%%%%%%%%%%%%%%%%%%%%%%%%%%%%%%%%%

\section{The geometric Demazure crystal and the Demazure quiver variety}
\label{sec:dqv}

In \cite{KS97,S02}, Kashiwara and Saito defined a crystal structure
on the set of irreducible components $B(\d)$ of $\sqcup_\v
\L(\v,\d)$, endowing it with the structure of the crystal of
$V(\lambda_\d)$.  As mentioned in Section~\ref{sec:dem-def}, there
is a Demazure crystal $B_w(\d)$ corresponding to the Demazure module
$V_w(\lambda_\d)$ for $w \in W$. Note that $B_w(\d)$ is isomorphic
to $B_w(\lambda_\d)$. For a fixed $\d$, let $\v_w$, $w \in W$ be
defined by
\[
\lambda_\d - \alpha_{\v_w} = w \lambda_\d.
\]

By Proposition~\ref{prop:qv-point}, $\L(\v_w,\d)$ is a point for all
$w \in W$.  Denote this point by $X_w = [x_w,t_w]$, the $G_\V$-orbit
of the representation $(x_w,t_w)$ of the quiver.

\begin{prop}
\label{prop:demazure-crystal}
For $X \in B(\d)$, we have that $X \in B_w(\d)$ if and only if all
points of $X$ are (orbits of) subrepresentations (up to
isomorphism) of $(x_w,t_w)$.
\end{prop}

\begin{proof}
We first prove that for all $X \in B_w(\d)$, $X$ consists of (orbits
of) subrepresentations (up to isomorphism) of $(x_w,t_w)$.  Our
proof is by induction on the length of $w$.   If this length is
zero, the statement is trivial.  Let $w=r_{i_l} \cdots r_{i_1}$ be a
reduced expression and $X \in B_w(\d)$.  If $X_{\lambda_\d} =
\mathfrak{L}(0,\d)$ then, by (1.2), we have
\[
    X=\tilde f_{i_l}^{n_l} \cdots \tilde f_{i_1}^{n_1} X_{\lambda_\d}
\]
for some $n_1, \dots, n_1 \in \Z_{\ge 0}$.  Let
\[
    X' = \tilde e_{i_l}^{n_l} X = \tilde f_{i_{l-1}}^{n_{l-1}}
    \cdots \tilde f_{i_1}^{n_1} X_{\lambda_\d}.
\]
Since $r_{i_l} w = r_{i_{l-1}} \cdots r_{i_1}$ is a reduced
expression, we see that $X' \in B_{r_{i_l}w}(\d)$ by (1.2).  It
follows from the inductive hypothesis that $X'$ consists of (orbits
of) subrepresentations (up to isomorphism) of $(x_{r_{i_l}w},
t_{r_{i_l}w})$.  Let
\[
    X_w = \{[x_w,t_w]\} \in B_w(\d),\quad X_{r_{i_l} w}
    = \{[x_{r_{i_l}w}, t_{r_{i_l}w}]\} \in B_{r_{i_l}w}(\d).
\]
Let $\V$, $\V'$, $\V^w$ and $\V^{r_{i_l} w}$ be the spaces
corresponding to the representations whose orbits are points of $X$,
$X'$, $X_w$ and $X_{r_{i_l} w}$ respectively.  By the definition of
the crystal operators in \cite{S02}, we have $\V^w_i = \V^{r_{i_l}
w}_i$ for $i \ne i_l$ and
\begin{equation} \label{eq:tildeVil}
    \V^w_{i_l} \xrightarrow[\cong]{((x_w)_h,t_w)} \ker \left(
    \bigoplus_{h\, :\, \inc(h)=i_l}
    \V^{r_{i_l} w}_{\out(h)} \oplus \W_{i_l} \xrightarrow{(\varepsilon(h)
    (x_{r_{i_l}w})_h,0)} \V^{r_{i_l} w}_{i_l} \right) \stackrel{\text{def}}{=}
    \tilde K.
\end{equation}

Since $X=\tilde f_{i_l}^{n_l} X'$, by the definition of the crystal
operators given in \cite{S02} we have
\[
    \V_i = \V'_i,\ i \ne i_l,\quad \text{and} \quad \V_{i_l} = \V'_{i_l}
    \oplus \C^{n_l},
\]
and an open dense subset of $X$ consists of orbits of
representations $(x,t)$ such that $\V'$ is $x$-stable and $[x',t']
\in X'$ where $x'=x|_{\V'}$, $t'=t|_{\V'}$. Now, by the stability
and moment map conditions, we have that the map
\begin{equation} \label{eq:Vil}
    \V_{i_l} \stackrel{(x_h,t)}{\hookrightarrow} \ker \left(
    \bigoplus_{h\, :\, \inc(h)=i_l} \V'_{\out(h)} \oplus
    \W_{i_l} \xrightarrow{(\varepsilon(h) x'_h,0)}
    \V'_{i_l} \right) \stackrel{\text{def}}{=} K
\end{equation}
is injective.  Since $(x',t')$ is a subrepresentation (up to
isomorphism) of $(x_{r_{i_l}w}, t_{r_{i_l}w})$, we have (after
replacing representations by different orbit representatives if
necessary) $K \subseteq \tilde K$. Thus, for $(x,t)$ in an open
dense subset of $X$, $(x,t)$ is isomorphic to a subrepresentation of
$(x_w,t_w)$ by \eqref{eq:tildeVil}, \eqref{eq:Vil} and the inductive
hypothesis.  Thus, to prove that $X$ consists entirely of orbits of
subrepresentations of $(x_w,t_w)$, it suffices to prove that the set
of such orbits is closed.  This set is precisely $\pi_1
(\pi_2^{-1}(\L(\v_w,\d))) = \pi_1 (\pi_2^{-1}([x_w,t_w]))$ where we
have used diagram \eqref{eq:diag_action_mod} with
\[
\v' = \v - \sum_{j=1}^l n_j \e^{i_j}.
\]
Since the map $\pi_2$ is proper (see \cite{N94,N98}) and the point
$[x_w,t_w]$ is obviously compact, we have that
$\pi_2^{-1}([x_w,t_w])$ is compact and thus $\pi_1
(\pi_2^{-1}([x_w,t_w]))$ is closed.

Now suppose that every point of $X$ is an orbit of a
subrepresentation of $(x_w,t_w)$ and that
\begin{equation}
\label{eq:Xw}
X_w = \kf_{i_l}^{m_l} \cdots \kf_{i_1}^{m_1} X_{\lambda_\d},
\end{equation}
Where $X_{\lambda_\d} = \L(0,\d)$.  If $l=0$ or, equivalently, $w
= \text{id}$ then the result follows.  So we assume $l > 0$.  From
the proof of Proposition~\ref{prop:qv-point}, we see that
\[
\ke_{i_l}^{m_l} X_w = X_{r_{i_l} w}.
\]
Note that $r_{i_l} w = r_{i_{l-1}} \dots r_{i_1}$ is a reduced
expression.

Let $\V^w$ be the space corresponding to the representation whose
orbit is $X_w$. We will use the notation $\V$ to denote the
subspace of $\V^w$ corresponding to a representation whose orbit
is a point of $X$. Note that this subspace depends on the point of
$X$ in question. From \eqref{eq:Xw} and the definition of the
crystal operators, we have an $x_w$-stable flag $0 \subset \V^1
\subset \dots \subset \V^l = \V^w$ such that $\dim \V^j / \V^{j-1}
= m_j \e^{i_j}$.

By the definition of the operator $\ke_{i_1}^{m_l}$ given in
\cite{S02}, $\ke_{i_l}^{m_l} X_w$ consists of the point $[x,t] =
[x_{r_{i_l}w},t_{r_{i_l}w}]$ where $(x,t)$ is the restriction of
$(x_w,t_w)$ to the space $\V^{l-1}$. We have $\V^{l-1}_i = \V^l_i
= \V^w_i$ for $i \ne i_l$ and
\[
\V^{l-1}_{i_l} = \im \left( \bigoplus_{\inc(h) = i_l}
\V^w_{\out(h)} \stackrel{x_w}{\longrightarrow} \V^w_{i_l} \right).
\]
Now, let $c$ be the maximum integer such that $\ke_{i_l}^c X \ne
0$.  By the definition of the operator $\ke_{i_1}^c$ given in
\cite{S02}, an open dense subset of $\ke_{i_l}^c X$ consists of
points $[x,t]$ where $(x,t)$ is the restriction of $(x',t')$, for
some $[x',t'] \in X$, to the space
\[
\im \left( \bigoplus_{\inc(h) = i_l} \V_{\out(h)}
\stackrel{x_w}{\longrightarrow} \V_{i_l} \right)
\]
at the $i_l$th vertex. Since $\V_{\out(h)} \subset \V^w_{out(h)}$
for any $h \in H$, the above is a subspace of $\V^{l-1}_{i_l}$.
Thus we have that (a representative of) every point of
$\ke_{i_l}^c X$ is a subrepresentation of (a representative of)
the unique point of $\ke_{i_l}^{m_l} X_w = X_{r_{i_l} w}$.  The result
then follows by induction on the length of $w$.
\end{proof}

\begin{defin}
For $w \in W$, we define $\L_w(\v,\d)$ to be the set of all $[x,t]
\in \L(\v,\d)$ such that $(x,t)$ is isomorphic to a
subrepresentation of $(x_w,t_w)$.  We call $\L_w(\v,\d)$ the
\emph{Demazure quiver variety}.
\end{defin}

%%%%%%%%%%%%%%%%%%%%%%%%%%%%%%%%%%%%%%%%%%%%%%%%%%%%%%%%%%%%%%%%%%%%%%%

\section{Demazure modules and the semicanonical basis}
\label{sec:sc-basis}

In \cite{Kas93}, Kashiwara proves that
\begin{equation}
\label{eq:global-basis}
V^q_w(\lambda) = \bigoplus_{b \in B_w(\lambda)} \Q(q)
G_\lambda(b),
\end{equation}
where $\{G_\lambda(b)\}_{b \in B(\lambda)}$ is the lower global
basis.  Thus the lower global basis behaves very nicely with
respect to the filtration of representations given by Demazure
modules.  It is known (cf. \cite{GL93}) that the lower global
basis coincides with Lusztig's canonical basis, defined
geometrically by considering perverse sheaves on quiver varieties.
Using the geometric description of the Demazure modules given
above, one could produce a geometric proof of
\eqref{eq:global-basis}.

In this section, we prove that the semicanonical basis also
behaves nicely with respect to the filtration by Demazure modules.
Let
\begin{equation}
L_w(\d) = \bigoplus_{X \in B_w(\d)} \C g_X.
\end{equation}

\begin{theo}
We have
\[
V_w(\lambda_\d) \cong L_w(\d).
\]
\end{theo}
\begin{proof}
We know that the Demazure module $V_w(\lambda_\d)$ is isomorphic
to the subset of constructible functions on $\bigsqcup_\v
\L(\v,\d)$ spanned by the result of acting on constant functions
on the point $X_w$ by all combinations of the $E_i$.  We also
denote this space of functions by $V_w(\lambda_\d)$.  Let $f$ be a
function in $V_w(\lambda_\d)$. By the definition of the action of
the $E_i$'s (see Section~\ref{sec:structure}), the support of $f$
consists of points $[x,t]$ where $(x,t)$ is a subrepresentation of
$(x_w,t_w)$.  Thus $f$ must be a linear combination of $g_X$ with
$X \in B_w(\d)$.  This is because for $X \not \in B_w(\d)$, $g_X$
is non-zero on an open dense subset of $X$ and thus by
Proposition~\ref{prop:demazure-crystal} contains a point $[x,t]$
in its support such that $(x,t)$ is not isomorphic to a
subrepresentation of $(x_w,t_w)$.  So we have that
$V_w(\lambda_\d) \subset L_w(\d)$.  However, by
\eqref{eq:demazure-crystalbase} or \eqref{eq:global-basis}, the
number of elements in the Demazure crystal $B_w(\d)$ is equal to
the dimension of the Demazure module.  Thus the result follows by
dimension considerations.
\end{proof}

So we see that the semicanonical basis has a compatibility with the
filtration by Demazure modules analogous to that of the canonical
basis. Note, however, that it was shown in \cite{GLS04} that the
semicanonical and (the specialization at $q=1$ of the) canonical
bases do indeed differ in general.

%%%%%%%%%%%%%%%%%%%%%%%%%%%%%%%%%%%%%%%%%%%%%%%%%%%%%%%%%%%%%%%%%%%%%%%

\section{Young pyramids and Demazure crystals for $\widehat{\mathfrak{sl}}_2$}
\label{sec:pyramids}

In \cite{Sav03b}, using the geometric realization of crystals
using quiver varieties, the author developed a combinatorial
realization of the crystals of integrable highest weight modules
of type $A_n^{(1)}$ using objects called Young pyramids.  These
can be considered a higher level analogue of the Young wall
realization of level one modules developed by Kang (cf.
\cite{K03}).  We note that Kang and Lee have also developed a
higher level generalization of Young walls, called Young slices
(cf. \cite{KL03}).  In this section, we consider the case $\g =
\widehat{\mathfrak{sl}}_2$ and give an explicit description of the
Demazure crystals in the language of Young pyramids.

Let $\alpha_0$ and $\alpha_1$ be the two simple roots of $\g$ and
let $r_i$ be the reflection with respect to $\alpha_i$.  The Weyl
group $W$ is generated by the $r_i$ and for every $n > 0$, $W$
contains two elements of length $n$.  They are
\begin{align*}
w_n^+ &= r_{i_n} \dots r_{i_2} r_{i_1},\quad i_j \equiv j+1 \mod
2, \\
w_n^- &= r_{i_n} \dots r_{i_2} r_{i_1},\quad i_j \equiv j \mod 2.
\end{align*}

\begin{figure}
\centering \epsfig{file=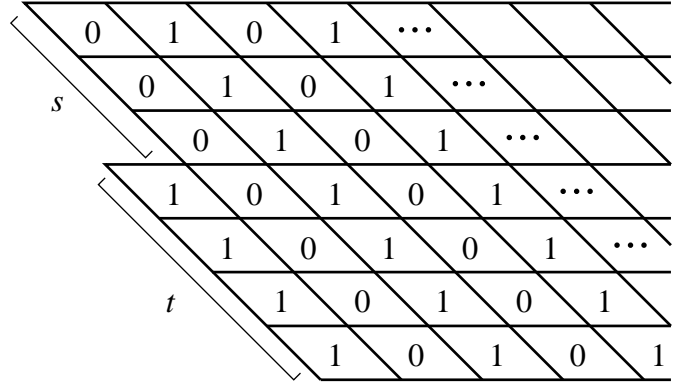,width=0.7\textwidth}
\caption{The ground state pyramid $P_\lambda$ of highest weight
$\lambda = s\omega_0 + t\omega_1$. \label{fig:groundpyramid}}
\end{figure}

Let $\lambda = s\omega_0 + t\omega_1$.  Then the ground state
pyramid $P_\lambda$ looks as in Figure~\ref{fig:groundpyramid}.
Let $\mathcal{P}(\lambda)$ be the set of all $1$-reduced Young
pyramids built on the ground state pyramid $P_\lambda$.  Recall
that Young pyramids are built by placing blocks of color $0$ or
$1$ on the ground state pyramid such that any block placed in an
empty slot matches the color of that slot and any block placed on
top of another block must be of opposite color.  The collection of
blocks placed in one slot is called a \emph{stack} and the words
\emph{row} and \emph{column} refer to positions in the ground
state pyramid.  We use the compass directions to refer to the
relative positions of slots in the ground state pyramid (or the
stacks placed on them). We number the columns, starting at zero,
from left to right. We adopt the convention that the westernmost
column is the column containing the $t$ slots of color 1. If $t=0$
then the westernmost column, containing $s$ slots of color 0 is
column one (and not zero).

\begin{figure}
\centering \epsfig{file=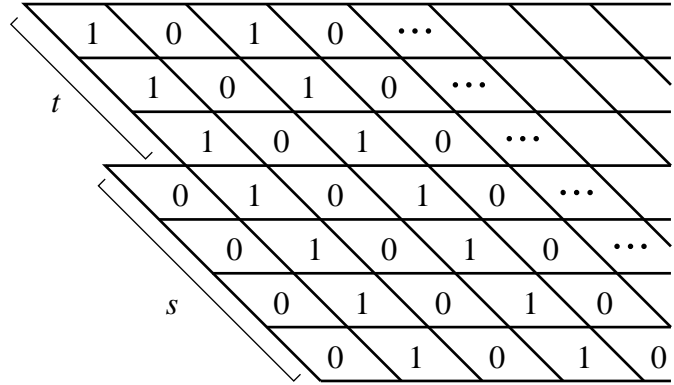,width=0.7\textwidth}
\caption{The modified ground state pyramid of highest weight
$\lambda = s\omega_0 + t\omega_1$ for describing the Demazure
crystal $\mathcal{P}_{w_n^+}$.
\label{fig:groundpyramid2}}
\end{figure}

Since the weight space of an extremal weight is one-dimensional,
we know that there is exactly one proper, $1$-reduced Young
pyramid in $\mathcal{P}(\lambda)$ of weight $w \lambda$ which we
will denote by $P_\lambda^w$. It is a straightforward calculation
based on weights to see that $P_\lambda^{w_n^-}$ is the Young
pyramid where the height of all stacks in the $i$th column is
$n-i$ for $0 \le i \le n$ and there are no other non-empty stacks.
For the case of $w_n^+$, we write the ground state pyramid for
$\lambda = s\omega_0 + t\omega_1$ as in
Figure~\ref{fig:groundpyramid2}.
The properness condition on Young pyramids ensures that the
crystal we obtain is the same with this modified ground state
pyramid. Then the unique pyramid $P_\lambda^{w_n^+}$ of weight
$w_n^+ \lambda$ is again the pyramid with stacks in the $i$th
column of height $n-i$ for $0 \le i \le n$ and no other stacks.
Let $\mathcal{P}_w(\lambda)$ be the subcrystal of
$\mathcal{P}(\lambda)$ corresponding to the Demazure crystal
$B_w(\lambda)$.  We also call $\mathcal{P}_w(\lambda)$ the
Demazure crystal.

\begin{theo}
The Demazure crystal $\mathcal{P}_{w_n^\pm}(\lambda)$ is the
subcrystal of $\mathcal{P}(\lambda)$ consisting of those Young
pyramids which are subpyramids of $P_\lambda^{w_n^\pm}$.
\end{theo}
\begin{proof}
We prove the case of $w_n^-$.  The case of $w_n^+$ is analogous. By
(1.2),
\[
    P_{w_n^-}(\lambda) = \{\kf_{i_n}^{k_n} \cdots \kf_{i_1}^{k_1} P_\lambda \ |\
    k_1, \dots, k_n \in \Z_{\ge 0},\ i_j \equiv j \mod 2 \}.
\]
Suppose $P \in P_{w_n^-}(\lambda)$. Since the blocks of any stack in
a Young pyramid must alternate color, we see that the maximum height
of any stack in $P$ is $n$ and the maximum height of any stack with
bottom block of color 0 (e.g. the stacks in the column one) is
$n-1$. Since $P$ is 1-reduced, the heights of stacks must strictly
decrease as we move east. Therefore, the height of any stack in the
$i$th column is less than or equal to $n-i$.  Therefore $P$ is a
subpyramid of $P_\lambda^{w_n^-}$.

Now let $P$ be a subpyramid of $P_\lambda^{w_n^-}$.  Then the
maximum possible height of a stack is $n$.  Suppose the maximum
height is $m$.  Since $P$ is $1$-reduced, the columns of height
$m$ must all be of the same color.  Let this color be $i$.  Recall
(see \cite{Sav03b}) that in defining the action of $\ke_i$ on a
Young pyramid, we compute the $i$-signature of $P$.  This is done
by arranging the stacks from tallest to shortest and assigning a
$-$ to $i$-removable stacks and a $+$ to $i$-admissible stacks. We
arrange these $+$ and $-$ signs from left to right (in the order
of the columns to which they correspond) and cancel $(+,-)$ pairs
until we obtain a sequence of $-$'s followed by $+$'s.  Then
$\ke_i$ acts by removing a block from the stack corresponding to
the rightmost $-$ sign (and $\ke_i P = 0$ if no such sign exists).
Now, the columns of height $m$, being of color $i$, cannot be
$i$-admissible.  Also, being the columns of maximal height, their
$-$ signs occur to the left of all the others and thus cannot be
cancelled in $(+,-)$ pairs.  Thus, after applying $\ke_i$ as many
times as possible (without obtaining zero), we see that the height
of all stacks is $\le m-1$.  We use here the fact that the only
way a stack of height $m$ can fail to be $i$-removable is if the
column immediately to the east is one block shorter.  Thus its top
block would also be of color $i$.  Continuing in this manner,
moving east, we see that we will eventually have an $i$-removable
block which cannot be cancelled in a $(+,-)$ pair.  When this
block is removed, the stack to the west becomes $i$-removable.
Continuing, the block on top of the stack of height $m$ will
eventually be removed after enough applications of the operator
$\ke_i$.

Now, the top of stacks which are now of the highest height must be
of color $(1-i)$ or else we could have applied $\ke_i$ more times.
So we repeat the above process, applying $\ke_{1-i}$ the maximum
number of times.  Since after each step we have reduced the
maximum height of the columns by at least one, we see that
\[
\ke_{i_k}^{c_k} \dots \ke_{i_1}^{c_1} P = P_\lambda,
\]
where $k \le n$ and the indices $i_1$ alternate between $0$'s and
$1$'s.  Therefore, by the description~\eqref{eq:dqv} of the
Demazure crystal, we see that $P \in \mathcal{P}_{w_n^-}(\lambda)$
which completes the proof.
\end{proof}

For $\lambda = s\omega_0 + t\omega_1$ (corresponding to $\d =
(\d_0,\d_1) = (s,t)$) and $w \in W$, let
\[
m_{\d,w} = \begin{cases} n & \text{if } w = w_n^-,\ t \ne 0 \\
n-1 & \text{if } w = w_n^-,\ t=0 \\
n & \text{if } w = w_n^+,\ s \ne 0 \\
n-1 & \text{if } w = w_n^+,\ s=0
\end{cases}.
\]

\begin{cor}
\label{cor:subpyramids}
$\mathcal{P}_w(\lambda_\d)$ is the subcrystal of
$\mathcal{P}(\lambda_\d)$ consisting of Young pyramids whose
stacks are all of height $\le m_{\d,w}$.
\end{cor}

Given the above descriptions of the Demazure crystals, the
character of the Demazure modules is given by
\begin{equation}
\dim V_{w}(\lambda) = \sum_{P \in \mathcal{P}_w(\lambda)}
e^{\wt(P)}.
\end{equation}
Explicit expressions for these characters were given in
\cite{FMO98} by characterizing the crystal bases in terms of
``paths". Certain specialized characters, obtained by applying the
map $e^\mu \mapsto q^{l(\mu)}$ where $l$ is some integral linear
function on the weight lattice, were obtained in \cite{San96}. The
dimensions of the Demazure modules for $\widehat{\mathfrak{sl}}_2$
were also computed in \cite{San96b}.  As an application of our
description in terms of Young pyramids, one can reproduce the
dimension formulas of \cite{San96b} by a simple counting argument.
Namely, we can show by counting the Young pyramids of
$\mathcal{P}_w(\lambda)$ that the dimensions of the Demazure
modules of $\widehat{\mathfrak{sl}}_2$ are given by
\begin{equation}
\dim V_w(s \omega_0 + t \omega_1) = \begin{cases}
(s+1)(s+t+1)^{n-1} & \text{if $w=w_n^+$ for some $n \ge 1$}, \\
(t+1)(s+t+1)^{n-1} & \text{if $w=w_n^-$ for some $n \ge 1$}.
\end{cases}
\end{equation}

%%%%%%%%%%%%%%%%%%%%%%%%%%%%%%%%%%%%%%%%%%%%%%%%%%%%%%%%%%%%%%%%%%%%%%%

\section{The Demazure quiver variety for
$\widehat{\mathfrak{sl}}_2$}
\label{sec:sl2-dqv}

In this section, we examine the Demazure quiver variety defined in
Section~\ref{sec:dqv} in the special case where $\g =
\widehat{\mathfrak{sl}}_2$.

For $\V \in \mathcal{V}$ and $x = (x_h)_{h \in H} \in
\mathbf{E}_\V$, we say that $x^n = 0$ if for any sequence $h_1,
\dots, h_n \in H$ such that $\inc(h_i) = \out(h_{i+1})$, we have
$x_{h_n} \dots x_{h_1} = 0$.

\begin{theo}
\label{thm:dqv-sl2}
When $\g = \widehat{\mathfrak{sl}}_2$,
\[
\L_w(\v,\d) = \{[x,t] \in \L(\v,\d)\ |\ x^{m_{\d,w}}=0\}.
\]
\end{theo}

Thus to obtain the Demazure quiver variety from the usual quiver
variety, we simply have to impose the condition that the quiver
representation is nilpotent of the given order.

\begin{proof}
Let $m=m_{\d,w}$.  It is easy to see from the description of the
extremal Young pyramid $P_w$ and the association of this pyramid
to the point of the quiver variety as described in \cite{Sav03b}
that $x_w^m = 0$.  Then the fact that all the points $[x,t]$ of
the Demazure quiver variety $\L_w(\v,\d)$ satisfy $x^m=0$ follows
from the fact that $(x,t)$ is a subrepresentation of $(x_w,t_w)$.

Now, let $[x,t] \in \L(\v,\d)$ with $x=(x_h)_{h \in H}$ and
$x^m=0$.  Then $x$ is in the conormal bundle to the $G_\V$-orbit
through $x_\Omega = (x_h)_{h \in \Omega}$, which corresponds to a
Young pyramid (see \cite{Sav03b}).  Since $x^m=0$, each stack in
the Young pyramid must have height at most $m$.  Now, by
Corollary~\ref{cor:subpyramids}, $\mathcal{P}_w(\lambda)$ is
precisely the set of Young pyramids with stacks of height $\le m$.
Thus, since the irreducible components of $\L(\v,\d)$ correspond
to the closures of the conormal bundles to orbits through the
$x_\Omega$ corresponding to Young pyramids, we see that $[x,t]$ is
a point of an irreducible component in $B_w(\d)$.  Then by
Proposition~\ref{prop:demazure-crystal}, we know that $(x,t)$ is
isomorphic to a subrepresentation of $(x_w,t_w)$ and therefore is
contained in $\L_w(\v,\d)$.
\end{proof}

\begin{cor}
Let $\text{irr}_X$ be the irreducible component of $\L_w(\v,\d)$
corresponding to $X \in B_w(\d)$.  Then for the Lie algebra
$\widehat{\mathfrak{sl}}_2$, we have
\[
\L_w(\v,\d) = \bigcup_{X \in B_w(\d)} \text{irr}_X.
\]
\end{cor}
\begin{proof}
This follows from the proof of Theorem~\ref{thm:dqv-sl2}.
\end{proof}

%%%%%%%%%%%%%%%%%%%%%%%%%%%%%%%%%%%%%%%%%%%%%%%%%%%%%%%%%%%%%%%%%%%%%%%

\bibliographystyle{abbrv}
\bibliography{biblist}

\end{document}